\documentclass{gtart}

\usepackage{amsmath,amsfonts}
\usepackage{pstricks,pst-node}

\def\ifplaintex{\expandafter\ifx\csname documentclass\endcsname\relax}

\expandafter\ifx\csname epsfbox\endcsname\relax\input epsf\fi

\ifplaintex 
\else
\fi

\def\gt{{\mathsurround=0pt\it $\cal G\mskip-2mu$eometry \&\ 
$\cal T\!\!$opology}}        

\def\gtp{{\mathsurround=0pt\it $\cal G\mskip-2mu$eometry \&\ 
$\cal T\!\!$opology $\cal P\!$ublications}}  

\def\volumenumber#1{\def\thevolumenumber{#1}}
\def\papernumber#1{\def\thepapernumber{#1}}
\def\volumeyear#1{\def\thevolumeyear{#1}}

\def\pagenumbers#1#2{\def\startpage{#1}\def\finishpage{#2}}
\def\published#1{\def\publishdate{#1}}
\def\proposed#1{\def\theproposer{#1}}
\def\seconded#1{\def\theseconders{#1}}
\def\received#1{\def\receiveddate{#1}}
\def\revised#1{\def\reviseddate{#1}}
\def\accepted#1{\def\accepteddate{#1}}

\long\def\asciiabstract#1{\long\def\theasciiabstract{#1}}

\let\\\par\let\thevolumenumber\relax\let\thepapernumber\relax
\let\thevolumeyear\relax\let\thesamplenumber\relax\let\startpage\relax
\let\finishpage\relax\let\publishdate\relax\let\receiveddate\relax
\let\reviseddate\relax\let\accepteddate\relax\let\theasciititle\relax
\let\theasciiauthors\relax
\let\theasciiabstract\relax
\let\theasciiemail\relax\let\theshortauthors\relax\let\theshorttitle\relax

\long\def\maketitlep{   

\count0=\startpage

\gt\hfill      
\hbox to 77pt{\vbox to 0pt{\vglue -15pt\epsfbox{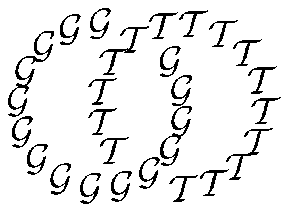}\vss}\hss}
\break
{\small\ifx\thesamplenumber\relax 
Volume \else Sample
\fi\thevolumenumber\ (\thevolumeyear)
\startpage--\finishpage\nl
Published: \publishdate}
\vglue 0.5truein plus 0.4fil minus 0.1truein

{\parskip=0pt\leftskip 0pt plus 1fil\def\\{\par\smallskip}{\ifplaintex\large
\else\Large\fi\bf\thetitle}\par\medskip}   

\vglue 0pt plus 0.1fil 

{\parskip=0pt\leftskip 0pt plus 1fil\def\\{\par}{\sc\theauthors}
\par\medskip}

\vglue 0pt plus 0.1fil 

{\small\parskip=0pt\let\newline\\
{\leftskip 0pt plus 1fil\def\\{\par}{\sl\theaddress}\par}
\expandafter\ifx\theemail\relax    
\relax\else\vglue 5pt plus 0.02fil minus 2pt\def\\{\stdspace{\rm 
and}\stdspace} 
\cl{Email:\stdspace\tt\theemail}\fi
\ifx\theurl\relax                  
\relax\else\vglue 5pt plus 0.02fil minus 2pt\def\\{\stdspace{\rm 
and}\stdspace}
\cl{URL:\stdspace\tt\theurl}\fi\par}

\vglue 7pt plus 0.3fil minus 3pt

{\bf Abstract}
\vglue 5pt plus 0.1fil minus 2pt

\theabstract

\vglue 7pt plus 0.3fil minus 3pt

{\bf AMS Classification numbers}\quad Primary:\quad \theprimaryclass

Secondary:\quad \thesecondaryclass

\vglue 5pt plus 0.3fil minus 2pt

{\bf Keywords:}\quad \thekeywords

\vglue 10pt plus 0.5fil minus 5pt

{\small  Proposed: \theproposer\hfill Received: \receiveddate\nl
Seconded: \theseconders\hfill 
\ifx\reviseddate\relax                         
Accepted: \accepteddate                        
\else
Revised: \reviseddate                          
\fi}
\eject
}       

\let\maketitlepage\maketitlep
\let\maketitle\maketitlepage

\font\phead=cmsl9 scaled 950
\font=cmsl9 scaled 1050
\font\pnum=cmbx10 scaled 913
\font\lnum=cmbx10 
\font\pfoot=cmsl9 scaled 950
\font=cmsl9 scaled 1050
\ifplaintex
\headline{\vbox to 0pt{\vskip -4.5mm\line{\small\phead\ifnum
\count0=\startpage ISSN 1364-0380 (on line)
1465-3060 (printed) \hfill {\pnum\folio}\else\ifodd\count0\def\\{ }%
\ifx\theshorttitle\relax\thetitle\else\theshorttitle\fi\hfill{\pnum\folio}
\else\def\\{ and }{\pnum\folio}\hfill\ifx\theshortauthors\relax\theauthors
\else\theshortauthors\fi\fi\fi}\vss}}
\footline{\vbox to 0pt{\vglue 0mm\line{\small\pfoot\ifnum\count0=\startpage
\copyright\ \gtp\hfill\else
\gt, Volume \thevolumenumber\ (\thevolumeyear)\hfill\fi}\vss
}}
\else
\makeatletter
\def\@oddhead{{\small\lhead\ifnum\count0=\startpage ISSN 1364-0380 (on line)
1465-3060 (printed) \hfill {\lnum\number\count0}\else\ifodd\count0
\def\\{ }\ifx\theshorttitle\relax \thetitle \else\theshorttitle\fi\hfill
{\lnum\number\count0}\else\def\\{ and }{\lnum\number\count0}
\hfill\ifx\theshortauthors\relax 
\theauthors\else\theshortauthors\fi\fi\fi}}\def\@evenhead{@oddhead}
\def\@oddfoot{\small\lfoot\ifnum\count0=\startpage\copyright\ \gtp\hfill\else
\gt, Volume \thevolumenumber\ (\thevolumeyear)\hfill\fi}
\def\@evenfoot{@oddfoot}
\makeatother
\fi

\newwrite\gtoutfile
\long\gdef\makeheadfile{  
{\def\\{, }
\immediate\openout\gtoutfile head.xxx
\immediate\write\gtoutfile{To: math@arxiv.org}
\immediate\write\gtoutfile{Subject: put}
\immediate\write\gtoutfile{--text follows this line--}
\immediate\write\gtoutfile{Proxy-for: \ifx\theasciiauthors\relax
\theauthors\else\theasciiauthors\fi <\ifx\theasciiemail\relax\theemail\else\theasciiemail\fi>}
\immediate\write\gtoutfile{\noexpand\\}
\immediate\write\gtoutfile{Authors: \ifx\theasciiauthors\relax
\theauthors\else\theasciiauthors\fi}
{\def\\{ }\immediate\write\gtoutfile{Title: \ifx\theasciititle\relax
\thetitle\else\theasciititle\fi}}
\immediate\write\gtoutfile{Subj-class: GT}
\immediate\write\gtoutfile{MSC-class: \theprimaryclass\ifx\thesecondaryclass\relax\else, \thesecondaryclass\fi}
\immediate\write\gtoutfile{Journal-ref: Geom. Topol. \thevolumenumber
(\thevolumeyear) \startpage-\finishpage}
\immediate\write\gtoutfile{Comments: Published in Geometry and Topology at}
\immediate\write\gtoutfile{    http://www.maths.warwick.ac.uk/gt/GTVol\thevolumenumber/paper\thepapernumber.abs.html}
\immediate\write\gtoutfile{\noexpand\\}
\immediate\write\gtoutfile{}
\ifx\theasciiabstract\relax
\immediate\write\gtoutfile{\theabstract}\else
\immediate\write\gtoutfile{\theasciiabstract}\fi
\immediate\write\gtoutfile{}
\immediate\write\gtoutfile{\noexpand\\}
\immediate\write\gtoutfile{}
\immediate\write\gtoutfile{<uuencoded .tar.gz file here>}
\immediate\write\gtoutfile{}
\immediate\closeout\gtoutfile}}  

\def\maketitlepage{\maketitlep\makeheadfile}
\let\maketitle\maketitlepage

\volumenumber{4}\papernumber{9}\volumeyear{2000}
\pagenumbers{277}{292}
\proposed{Robion Kirby}
\seconded{Michael Freedman, Walter Neumann}
\received{4 August 1999}
\revised{28 September 2000}
\accepted{21 September 2000}
\published{8 October 2000}

\newtheorem{theorem}{Theorem}
\newtheorem{lemma}[theorem]{Lemma}
\newtheorem{corollary}[theorem]{Corollary}
\newtheorem{conjecture}[theorem]{Conjecture}
\newtheorem{question}[theorem]{Question}
\theoremstyle{remark}
\newtheorem{example}[theorem]{Example}

\let\citealp\cite

\newcommand{\R}{{\mathbb{R}}}
\newcommand{\Z}{{\mathbb{Z}}}
\renewcommand{\H}{{\mathbb{H}}}
\newcommand{\cP}{{\mathcal{P}}}
\newcommand{\cQ}{{\mathcal{Q}}}
\newcommand{\cA}{{\mathcal{A}}}
\newcommand{\cB}{{\mathcal{B}}}
\newcommand{\cC}{{\mathcal{C}}}
\newcommand{\se}{\succeq}
\newcommand{\pe}{\preceq}
\newcommand{\Ie}{Ie}
\newcommand{\eg}{eg}
\newcommand{\ie}{ie}
\DeclareMathOperator{\Vol}{Vol\ }

\newenvironment{fullfigure}[2]
    {\begin{figure}[ht!]\small\begin{center}\def\ffa{#1}\def\ffb{#2}}
    {\vspace{\baselineskip}\caption{\ffb}\label{\ffa}\end{center}\end{figure}}

\newcommand{\swslug}
    {\pscustom[fillcolor=gray15,fillstyle=solid]
    {\psline(-.75,-.75)(-.75,.75)(.75,.75)(.75,-.75)(-.75,-.75)
    \psline[liftpen=2](-.20,-.60)
    (.60,-.60)(.60,.60)(-.60,.60)(-.60,-.20)
    \psarcn(-.6,-.6){.4}{90}{0}}}
\newcommand{\seslug}{\rput{90}(0,0){\swslug}}
\newcommand{\neslug}{\rput{180}(0,0){\swslug}}
\newcommand{\nwslug}{\rput{270}(0,0){\swslug}}

\psset{linewidth=.5pt,unit=.25in,arrowsize=2pt 3}
\newgray{gray15}{.85}
\SpecialCoor

\begin{document}
\title{Notions of denseness}
\author{Greg Kuperberg}
\address{Department of Mathematics, University of California\\One
Shields Ave, Davis, CA 95616-8633, USA}
\email{greg@math.ucdavis.edu}
\begin{abstract}
The notion of a completely saturated packing \cite{Kuperberg:saturated} is a
sharper version of maximum density, and the analogous notion of a completely
reduced covering is a sharper version of minimum density. We define two related
notions:  uniformly recurrent and weakly recurrent dense packings, and
diffusively dominant packings.  Every compact domain in Euclidean space has a
uniformly recurrent dense packing.  If the domain self-nests, such a packing is
limit-equivalent to a completely saturated one. Diffusive dominance is yet
sharper than complete saturation and leads to a better understanding of
$n$--saturation.
\end{abstract}
\asciiabstract{The notion of a completely saturated packing [Fejes Toth, 
Kuperberg and Kuperberg, Highly saturated packings and reduced
coverings, Monats. Math. 125 (1998) 127-145] is a sharper version of
maximum density, and the analogous notion of a completely reduced
covering is a sharper version of minimum density. We define two
related notions: uniformly recurrent and weakly recurrent dense
packings, and diffusively dominant packings.  Every compact domain in
Euclidean space has a uniformly recurrent dense packing.  If the
domain self-nests, such a packing is limit-equivalent to a completely
saturated one. Diffusive dominance is yet sharper than complete
saturation and leads to a better understanding of n-saturation.}

\keywords{Density, saturation, packing, covering, dominance}

\primaryclass{52C15, 52C17}
\secondaryclass{52C20, 52C22, 52C26, 52B99}

\maketitle

What is the best way to pack congruent copies of some compact domain in a
geometric space?  Is there a best way?  The first notion of optimality is
\emph{maximum density}.  While this notion is useful, it is also incomplete. 
For example, the packing $\cA$ in Figure~\ref{f:circle} has the maximum
density, but there is clearly something missing.  Intuitively, the packing
$\cB$ is uniquely optimal.

\begin{fullfigure}{f:circle}
    {Two dense circle packings}
\makebox[1.5em]{$\cA$} \pspicture*[](-6,-3.464)(6,3.464)
\multirput(-6,3.464)(1,0){13}{\pscircle[fillcolor=gray15,fillstyle=solid](0,0){.5}}
\multirput(-5.5,2.598)(1,0){12}{\pscircle[fillcolor=gray15,fillstyle=solid](0,0){.5}}
\multirput(-6,1.732)(1,0){13}{\pscircle[fillcolor=gray15,fillstyle=solid](0,0){.5}}
\multirput(-5.5,.866)(1,0){5}{\pscircle[fillcolor=gray15,fillstyle=solid](0,0){.5}}
\multirput(1.5,.866)(1,0){5}{\pscircle[fillcolor=gray15,fillstyle=solid](0,0){.5}}
\multirput(-6,0)(1,0){5}{\pscircle[fillcolor=gray15,fillstyle=solid](0,0){.5}}
\pscircle[fillcolor=gray15,fillstyle=solid](0,0){.5}
\multirput(2,0)(1,0){5}{\pscircle[fillcolor=gray15,fillstyle=solid](0,0){.5}}
\multirput(-5.5,-.866)(1,0){5}{\pscircle[fillcolor=gray15,fillstyle=solid](0,0){.5}}
\multirput(1.5,-.866)(1,0){5}{\pscircle[fillcolor=gray15,fillstyle=solid](0,0){.5}}
\multirput(-6,-1.732)(1,0){13}{\pscircle[fillcolor=gray15,fillstyle=solid](0,0){.5}}
\multirput(-5.5,-2.598)(1,0){12}{\pscircle[fillcolor=gray15,fillstyle=solid](0,0){.5}}
\multirput(-6,-3.464)(1,0){13}{\pscircle[fillcolor=gray15,fillstyle=solid](0,0){.5}}
\psframe(-6,-3.464)(6,3.464)
\endpspicture \vspace{\baselineskip} \\
\makebox[1.5em]{$\cB$} \pspicture*[](-6,-3.464)(6,3.464)
\multirput(-6,3.464)(1,0){13}{\pscircle[fillcolor=gray15,fillstyle=solid](0,0){.5}}
\multirput(-5.5,2.598)(1,0){12}{\pscircle[fillcolor=gray15,fillstyle=solid](0,0){.5}}
\multirput(-6,1.732)(1,0){13}{\pscircle[fillcolor=gray15,fillstyle=solid](0,0){.5}}
\multirput(-5.5,.866)(1,0){12}{\pscircle[fillcolor=gray15,fillstyle=solid](0,0){.5}}
\multirput(-6,0)(1,0){13}{\pscircle[fillcolor=gray15,fillstyle=solid](0,0){.5}}
\multirput(-5.5,-.866)(1,0){12}{\pscircle[fillcolor=gray15,fillstyle=solid](0,0){.5}}
\multirput(-6,-1.732)(1,0){13}{\pscircle[fillcolor=gray15,fillstyle=solid](0,0){.5}}
\multirput(-5.5,-2.598)(1,0){12}{\pscircle[fillcolor=gray15,fillstyle=solid](0,0){.5}}
\multirput(-6,-3.464)(1,0){13}{\pscircle[fillcolor=gray15,fillstyle=solid](0,0){.5}}
\psframe(-6,-3.464)(6,3.464)
\endpspicture
\end{fullfigure}

The problem of defining optimality for packings and coverings has been
considered previously by Fejes T\'oth , Conway and Sloane  and probably by
others.   Fejes T\'oth \cite{Fejes-Toth:solid} defines a packing (or covering)
to be \emph{solid} if no finite rearrangement of its elements yields a
non-congruent packing (or covering).  This is a very strict condition which is
only sometimes achievable.  Conway and Sloane \cite{CS:best} propose a
definition of optimality for sphere packings only which, as we will discuss in
Section~\ref{s:others}, is too strict to be useful.  In this article we will
introduce two other notions:  uniformly and weakly recurrent dense packings, and
diffusively dominant packings.  We will show that these notions have some
favorable properties, and in particular we will compare them to the notion of
complete saturation \cite{Kuperberg:saturated}.

A packing $\cP$ of some body (\ie, a compact domain) $K$ is
\emph{$n$--saturated} if it is not possible to replace any $n-1$ elements of
$\cP$ by $n$ and still have a packing. It is
\emph{completely saturated} if it is saturated for all $n$.  We can express
saturation by a weak partial ordering on packings: $\cP_1 \se_{S,n} \cP_2$ if
the packing $\cP_1$ is obtained from $\cP_2$  by repeatedly replacing $n-1$
elements of $\cP_2$ by $n$ or more (either a larger finite number or an
infinite number).  Dropping the restriction on $n$, we say that $\cP_1 \se_S
\cP_2$. Thus $\cP$ is completely saturated if and only if it is maximal with
respect to the partial ordering $\se_S$.  For example, the packing $\cA$ in
Figure~\ref{f:circle} is not 1--saturated.  On the other hand, any packing
which is both periodic and has maximum density, such as the packing $\cB$
in Figure~\ref{f:circle}, is completely saturated.

Dually, a covering $\cP$ is \emph{$n$--reduced} if it is not possible to replace
$n$ elements by $n-1$ and still have a covering. It is \emph{completely
reduced} if it is $n$--reduced for all $n$.  The main open problem
about complete saturation and complete reduction is the following:

\begin{question} Does every body $K$ in Euclidean space $\R^d$ admit a
completely saturated packing and a completely reduced covering? What about
bodies in hyperbolic space $\H^d$ or other homogeneous spaces?
\end{question}

In this article we will present a stronger version of a previous result
\cite{Kuperberg:saturated} that strictly self-nesting domains have completely saturated
packings and completely reduced coverings.  By Theorem~\ref{th:saturated},
every dense packing of a self-nesting $K$ has arbitrarily large regions that
are arbitrarily close to some completely saturated packing.  The technical
argument is actually the same; only the conclusion is different.

A packing $\cA$ is a \emph{limit} of a packing $\cB$, or $\cA \se_L \cB$, if a
sequence of translates of $\cB$ converges to $\cA$ in the Hausdorff topology. 
It is \emph{uniformly recurrent} if it is maximal with respect to this weak
partial ordering in the space of packings.  It is \emph{weakly recurrent} if it
is a limit of (generally different) uniformly recurrent packings. Given a
finite measure $\mu$ in $\R^d$, a packing $\cA$ \emph{$\mu$--dominates} a
packing $\cB$, or $\cA \se_\mu \cB$, if the inequality
$$\mu * \chi_\cA \ge \mu * \chi_\cB$$
holds everywhere, where $\chi_\cP$ is the characteristic measure of $\cP$ and
$\mu * \nu$ is the convolution of the measures $\mu$ and $\nu$. We say that
$\cA$ \emph{diffusively dominates} $\cB$, or $\cA \se_D \cB$, if $\cA \se_\mu
\cB$ for some $\mu$.  Diffusive domination is also a weak partial
ordering, since 
$$\cA \se_\mu \cB \se_\nu \cC \implies \cA \se_{\mu * \nu} \cC.$$
The packing $\cA$ is \emph{$\mu$--dominant} or
\emph{diffusively dominant} if it is maximal with respect to the corresponding
partial orderings.  Diffusive dominance implies $\mu$--dominance
for all $\mu$ as follows:  If $\cA$ were the former but not the latter, then
$$\mu * \chi_\cA \le \mu * \chi_\cB$$
for some $\mu$ without equality holding everywhere, which would preclude
the possibility that $\cA \se_{\mu * \nu} \cB$ for any $\nu$.
Yet there would be a $\nu$ such that $\cA \se_\nu \cB$.

Unfortunately there exists a sequence of completely saturated packings of a
square that converges to a packing which is not even simply saturated
\cite{Kuperberg:saturated}.   In the present context we can conclude that the
saturation partial orderings $\se_{S,n}$ and $\se_S$ are irreparably
incompatible with the Hausdorff topology on the space of packings.  By
contrast, limits and uniform recurrence \emph{are} compatible, in the sense
that the set of pairs of packings $\cA \se_L \cB$ is closed. The $\mu$--dominant
partial orderings have an intermediate property: For any fixed $\mu$, the set
of pairs $\cA \se_\mu \cB$ is closed. Nonetheless the set of pairs $\cA \se_D
\cC$ is not closed.

Our arguments below apply to the dual case of thin coverings. Indeed, they
apply more generally to $k$--packings or $k$--coverings of finite collections of
bodies (or even compact families of bodies), which could have weighted density
or other variations. They also apply to restricted congruences, such as
packings by translates, provided that the set of allowed isometries includes
all translations.  Theorem~\ref{th:recurrent} also applies to packings and
coverings in arbitrary homogeneous Riemannian geometries, such as hyperbolic
geometry and complex hyperbolic geometry, while Theorem~\ref{th:constant}
applies to packings and coverings in sub-exponential geometries such as
nil-geometry and solve-geometry \cite{Thurston:geometry}.  We will assume these
generalizations implicitly.

\rk{Acknowledgments}
The author would like to thank W{\l}odzimierz Kuperberg and Bill Thurston for
helpful discussions pertaining to this article, and the referees
for useful remarks.

The author is supported by NSF grant DMS \#9704125 and by a Sloan
Foundation Research Fellowship.

\section{Basic notions}

A \emph{domain} is a set which is the closure of its interior and a \emph{body}
is a compact domain. Let $\cP$ be a collection of congruent copies of a body
$K$ in Euclidean space $\R^d$. The \emph{density} $\delta(\cP)$ of $\cP$ is
defined as 
\begin{equation}
\delta(\cP) = \lim_{r \to \infty}
\frac{\sum_{D \in \cP} \Vol B(p,r) \cap D}{\Vol B(p,r)}.
\label{e:density} \end{equation}
In this expression $p$ is an arbitrary (fixed) point in $\R^d$ and $B(p,r)$ is
the round ball of radius $r$ centered at $p$.  The collection $\cP$ has
\emph{uniform density} if the convergence in equation~\eqref{e:density} is
uniform in $p$.

The collection $\cP$ is a \emph{packing} if its elements have disjoint
interiors, and that it is a \emph{covering} if the union of its elements is all
of $\R^d$.  A packing is \emph{dense} if it has maximum density, and a covering
is \emph{thin} if it has minimum density. \emph{Uniformly dense} means that
density is both maximum and uniform.

These definitions have the following desirable properties
\cite{Groemer:existence}:
\begin{enumerate}
\item If the limit exists for some $p$, it exists for all $p$. 
\item The supremum $\delta(K)$ (the packing density of $K$) of the density
$\delta(\cP)$ is achieved by a uniformly dense packing.
\item $\delta(K)$ is greater than or equal to the upper density of any packing, in
which the limit in equation~\eqref{e:density} is replaced by a lim sup.   
\item If $\cP$ has uniform density, then
$$\delta(\cP) = \lim_{r \to \infty}
\frac{\sum_{D \in \cP} \Vol (rR+p) \cap D }{\Vol rR}$$
uniformly in $p \in \R^d$, where $R$ is an arbitrary body and $rR+p$ denotes
$R$ dilated by a factor of $r$ and translated by $p$.
\end{enumerate}

We define a modified Hausdorff distance between collections of domains and we
will use the induced Hausdorff topology.  By definition, $d(\cP,\cQ) \le
\epsilon$ means that there exists a binary relation $K \sim L$ for $K \in \cP$
and $L \in \cQ$ such that $K \sim L$ implies that $d(K,L) \le \epsilon$.  For
each $K \in \cP$ in the ball $B(0,1/\epsilon)$, there must exist a unique $L
\in \cQ$ such that $K \sim L$  and vice versa.  (The relation need not satisfy
any other axioms.)

A \emph{weak partial ordering} $\se$ on a set $S$ is defined as any binary
relation such that $x \se x$ and $x \se y \se z$ implies $x \se z$. Any weak
partial ordering induces an equivalence relation: $x$ and $y$ are equivalent
when $x \se y \se x$.

\section{Recurrence}
\label{s:recur}

All of the results of this section have implications for 
weakly recurrent packings, but it is more natural to state them
in terms of uniformly recurrent packings.

\begin{theorem} Every packing $\cP$ of a body $K$ admits a uniformly recurrent
packing $\cQ$ as a limit.  If $\cP$ is uniformly dense, so is $\cQ$.
\label{th:recurrent}
\end{theorem}
\begin{proof}
The translation group $\R^d$ acts on the space of packings in $\R^d$, and the
action is continuous relative to our Hausdorff topology.  This group
action may be interpreted as a dynamical system with a $d$--dimensional time
variable.  The notions of limits and uniform recurrence are then familiar from
the dynamical systems point of view.

To say that $\cA$ is a limit of $\cB$ is to say that $\cA$ lies in the closure
of the orbit of $\cB$. A limit equivalence class is an invariant set (not
necessarily closed) in which each orbit is dense.  If such an invariant set is
closed, it is called a \emph{minimal set}. A uniformly recurrent packing is a
point in a minimal set in this dynamical system.  Since the space of packings
is compact, minimal sets exists in the closure of every orbit.

The set of uniformly dense packings is unfortunately not compact, but it is
exhausted by invariant compacta in the following way:  If $\cP$ is uniformly
dense, the proportion of a large disk covered by $\cP$ must approach the
packing density $\delta(K)$ at a certain rate.  \Ie, there is some function
$\epsilon(r)$ which converges to 0 as $r$ goes to infinity and such that the
density in any ball $B(p,r)$ is within $\epsilon(r)$ of $\delta(K)$.  The set
of packings that are uniform using any fixed $\epsilon(r)$ is compact, and it
is also invariant under the action.
\end{proof}

Note that the density of a uniformly recurrent packing $\cP$ necessarily exists
uniformly.  In particular if $\cP$ is dense, it is uniformly dense.

Theorem~\ref{th:recurrent} can also be argued from the point of view of general
topology:  If $X$ is a compact Hausdorff space, and if $\se$ is a weak partial
ordering which is closed as a subset of $X \times X$, then it has a maximal
element.  (Proof: By compactness every ascending chain $C$ has a convergent
cofinal subchain.  By closedness, the limit is an upper bound for $C$.  
Therefore Zorn's Lemma produces a maximal element of $\se$ on $X$.)

Given a body $K$ in Euclidean space and an allowed isometry group of $G$ for
collections of copies $K$ (\eg, the group of translations or the group of all
isometries), the domain $K$ (strictly) \emph{self-nests} if for every $\epsilon
> 0$ there exists $g \in G$ such that $gK$ lies in the interior of
$(1+\epsilon)K$.  For collections of translates, the condition holds if $K$ is
strictly star-shaped (star-shaped with a continuous radial function).  For
general arrangements of congruent copies, it is a restatement of the strict
nested similarity property of Reference \citealp{Kuperberg:saturated}. For
example, the spiral-like body in Figure~\ref{f:spiral} (the 1--neighborhood of a
part of a logarithmic spiral) has this property.

\begin{fullfigure}{f:spiral}
    {A self-nesting spiral-like body}
\pspicture(-4,-3.7)(5,2.5)
\pscurve[arrows=c-c,linewidth=1.06]
    (0.15;-330)(0.16;-300)(0.18;-270)(0.20;-240)(0.22;-210)(0.24;-180)
    (0.26;-150)(0.29;-120)(0.32;-90)(0.35;-60)(0.39;-30)(0.42;0)
    (0.47;30)(0.51;60)(0.56;90)(0.62;120)(0.68;150)(0.75;180)
    (0.83;210)(0.91;240)(1.00;270)(1.10;300)(1.21;330)(1.33;360)
    (1.46;390)(1.61;420)(1.77;450)(1.95;480)(2.14;510)(2.36;540)
    (2.59;570)(2.85;600)(3.14;630)(3.45;660)(3.80;690)(3.92;700)(4.04;710)(4.18;720)
\pscurve[arrows=c-c,linewidth=1,linecolor=gray15]
    (0.15;-330)(0.16;-300)(0.18;-270)(0.20;-240)(0.22;-210)(0.24;-180)
    (0.26;-150)(0.29;-120)(0.32;-90)(0.35;-60)(0.39;-30)(0.42;0)
    (0.47;30)(0.51;60)(0.56;90)(0.62;120)(0.68;150)(0.75;180)
    (0.83;210)(0.91;240)(1.00;270)(1.10;300)(1.21;330)(1.33;360)
    (1.46;390)(1.61;420)(1.77;450)(1.95;480)(2.14;510)(2.36;540)
    (2.59;570)(2.85;600)(3.14;630)(3.45;660)(3.80;690)(3.92;700)(4.04;710)(4.18;720)
\pscurve
    (0.15;-330)(0.16;-300)(0.18;-270)(0.20;-240)(0.22;-210)(0.24;-180)
    (0.26;-150)(0.29;-120)(0.32;-90)(0.35;-60)(0.39;-30)(0.42;0)
    (0.47;30)(0.51;60)(0.56;90)(0.62;120)(0.68;150)(0.75;180)
    (0.83;210)(0.91;240)(1.00;270)(1.10;300)(1.21;330)(1.33;360)
    (1.46;390)(1.61;420)(1.77;450)(1.95;480)(2.14;510)(2.36;540)
    (2.59;570)(2.85;600)(3.14;630)(3.45;660)(3.80;690)(3.92;700)(4.04;710)(4.18;720)
\endpspicture
\end{fullfigure}

\begin{theorem} If $K \subset \R^d$ is a self-nesting body, then every
uniformly recurrent, dense packing $\cP$ is limit-equivalent to a completely
saturated one.
\label{th:saturated}
\end{theorem}

We follow the argument in Reference~\citealp{Kuperberg:saturated}.

\begin{lemma} Let $\cP$ be a dense packing of a self-nesting body $K \subset
\R^d$. For every radius $r$ and every $\epsilon > 0$, there is an $\alpha > 0$
with the following property:  Consider points $p$ such that the restriction of
the packing $\cP$ to the ball $B(p,r)$ is within $\alpha$ in modified Hausdorff
distance of a packing in this disk that is not completely saturated.  The set
$S$ of such points has density at most $\epsilon$. \label{l:technical}
\end{lemma}

\begin{proof} Since $K$ is self-nesting, we can \emph{loosen} any packing $\cP$
by homothetic expansion: We expand the packing and each copy of $K$ by a factor
of $1+\gamma$, and then we replace each copy of $(1+\gamma)K$ by a copy of  $K$
contained in its interior. For every $\gamma$ and $K$, there is an $\alpha$
such that in the loosened packing $\cP_\gamma$, no two copies of $K$ are within
$2\alpha$ of each other.

Intuitively, if we loosen $\cP$, we decrease the density by expansion, but we
may then increase it by re-saturation.  We choose the constant $\gamma$
so that the latter would outweigh the former if $\cP$ failed to satisfy the
conclusion of the theorem.  More precisely, we choose $\gamma$ so that
$$d\gamma \delta(K) < \frac{\epsilon (\Vol K)}{\Vol B(p,2r)}.$$ 
Let $\alpha$ be the constant in the previous paragraph. Consider a maximal
packing of balls $\{B(p_i,r)\}$ such that $\cP$ is within $\alpha$ of an
unsaturated packing in each ball.  By maximality, the corresponding collection
$\{B(p_i,2r)\}$ covers $S$.  In the loosened packing $\cP_\gamma$, we can cram
in at least one extra copy of $K$ in each expanded ball $(1+\gamma)B(p_i,r)$.
If $S$ had upper density greater than $\epsilon$, the new packing $\cP'_\gamma$
would have greater density than that of $\cP$, contradicting the assumption that $\cP$
is dense.
\end{proof}

\begin{proof}[Proof of Theorem~\ref{th:saturated}] Let $\epsilon_k =
1/2^{k+1}$ and $r_k = 2^k$, and choose the corresponding $\alpha_k$ according
to Lemma~\ref{l:technical}.  For each $n$, there is a non-zero measure of
points $p \in \R^d$ such that $\cP$ is simultaneously at least $\alpha_k$ away
from unsaturated in $B(p,r_k)$ for all $1 \le k \le n$. Let $p_n$ be one such
point.  The sequence of translated packings $\{\cP - p_n\}$ must have  a
convergent subsequence.  The limit $\cQ$ of this subsequence has the same
property for all $k$; in particular it is completely saturated.

By construction $\cQ$ is a limit of $\cP$. Since $\cP$ is uniformly recurrent,
the two packings are limit-equivalent.
\end{proof}

\begin{question} Is it possible that in an equivalence class of uniformly
recurrent packings of a self-nesting domain $K$, some representatives are
completely saturated and others are not?
\end{question}

\begin{theorem} The only uniformly recurrent dense packing of a circle
(a round disk) is the hexagonal packing.
\label{th:hex}
\end{theorem}
\begin{proof} The result follows from the fact that the unique smallest
possible Vor\-on\-oi region in a circle packing is a regular hexagon
\cite{Fejes-Toth:kugel}. For every $\epsilon > 0$, the union of the Voronoi
regions that are more than $\epsilon$ away from the regular circumscribed
hexagon must have density 0 in the plane in a dense packing.  Consequently
there are arbitrarily large disks in which every Voronoi region is within
$\epsilon$ of the optimal shape.  Taking the limit $\epsilon \to 0$, the
conclusion is that there are arbitrarily large regions that converge to the
hexagonal packing. Thus the hexagonal packing is a limit of every dense
packing.

Conversely, any periodic packing is the only limit of itself. In
particular, the hexagonal circle packing has this property.
\end{proof}

Note that even if there are distinct uniformly dense packings they may all form
a single minimal set, or limit equivalence class.  In this case the equivalence
class as a whole may be considered the unique uniformly recurrent solution to
the packing problem.

\begin{example} The Penrose tilings of unit rhombuses form a limit equivalence
class.  Non-Penrose tilings of the rhombuses are usually eliminated by certain
matching rules, but one may equally well use notches and teeth for this
purpose.  The notches and teeth may be chosen so that the tiles are
self-nesting. Indeed there is a set of three convex polygons due to Robert
Ammann that only tile in Penrose fashion \cite{GS:patterns}. Thus
Theorem~\ref{th:saturated} applies to families of bodies whose densest packings
fall into uncountable limit equivalence classes.
\end{example}

\begin{example} Bi-infinite sequences of symbols from a finite alphabet are
equivalent to ``packings'' of unit intervals colored by the same alphabet.
Every periodic sequence is uniformly recurrent, but there are also others, such
as the Thue--Morse sequence \cite{MH:unending} in the two letter alphabet
$\{0,1\} = \Z/2$. The $n$th term of the sequence, for $n \ge 0$, is the sum in
$\Z/2$ of the binary digits of $n$.  The sequence is only infinite in one
direction, but any bi-infinite limit of translates of the Thue--Morse
sequence is also
uniformly recurrent and aperiodic. By contrast, every sequence over a finite
alphabet is weakly recurrent, since every finite sequence can be extended to a
periodic one.  Binary sequences are also a model for Barlow packings of spheres
in $\R^3$ \cite{Barlow:probable}.  It seems likely that a sphere packing is
weakly recurrent among dense packings if and only if it is Barlow
\cite{Hales:overview,Sloane:confirmed}, but not all Barlow packings are
uniformly recurrent.
\end{example}

\section{Diffusion}
\label{s:diffusion}

We first state some of the notions of diffusion and measure more precisely. We
assume that all measures are defined on Borel sets in $\R^d$.  If a measure
$\mu$ has a density function $f(x)$,
$$\mu(A) = \int_A f(x) dx,$$
we will write $\mu(dx)$ for the density $f(x) dx$.  The characteristic measure
$\chi_S$ of a set $S$ is defined by
$$\chi_S(A) = \Vol A \cap S.$$
If $\cP$ is a collection of sets, its characteristic measure
$\chi_\cP$ is defined as the sum:
$$\chi_\cP = \sum_{A \in \cP} \chi_A.$$
If $\mu$ and $\nu$ are measures in $\R^d$ and the total measure of $\mu$ is
finite, then their convolution $\mu * \nu$ is  defined by linear extension of
addition of points:
$$(\mu * \nu)(A) = (\mu \times \nu)(\{(p,q) | p+q \in A\}),$$
where $\mu \times \nu$ is the product measure on $\R^{2d}$.

\begin{fullfigure}{f:string}
    {Anomalously replacing one square by two}
\pspicture(-2,.5)(8,6)
\rput(-1,5.5){$\cP_2$}\multirput(0,5)(1,0){2}
    {\psframe[fillstyle=solid,fillcolor=gray15](0,0)(1,1)}
\rput(-1,4){$\cP_3$}\multirput(1,3.5)(1,0){3}
    {\psframe[fillstyle=solid,fillcolor=gray15](0,0)(1,1)}
\rput(-1,2.5){$\cP_4$}\multirput(2,2)(1,0){4}
    {\psframe[fillstyle=solid,fillcolor=gray15](0,0)(1,1)}
\rput(-1,1){$\cP_5$}\multirput(3,.5)(1,0){5}
    {\psframe[fillstyle=solid,fillcolor=gray15](0,0)(1,1)}
\endpspicture
\end{fullfigure}

Even though the saturation partial orderings $\se_{S,n}$ are not closed in the
Hausdorff topology, we can still ask what happens if we have an increasing
sequence of packings. The idea is that the packings in the sequence are
improving, so perhaps they must converge to a packing which is better still. 
Unfortunately this is not the case.  For example, we can let $\cP_k$ be a
string of $k$ unit squares whose centers start at $(k,0)$ and ending at
$(2k-1,0)$ (Figure~\ref{f:string}). Although $\cP_{k+1}$ is obtained from
$\cP_k$ by replacing one square by two, the limit is the empty packing.

This disease may be cured by considering another partial ordering $\se_{C,n}$
on packings of a body $K$. We say that $\cA \se_{C,n} \cB$ if $\cA$ can be
obtained from $\cB$ by replacing $k<n$ copies of the body $K$ by $k+1$ in such
a way that the union of the $2k+1$ copies of $K$ is a connected set. Clearly
$\se_{C,n}$ and $\se_{S,n}$ have the same maximal elements.

\begin{theorem} If a sequence of packings $\{\cP_i\}$ of a body $K$ in
Euclidean space $\R^d$ increases under the connected $n$--saturation partial
ordering, $$\cP_1 \pe_{C,n} \cP_2 \pe_{C,n} \cP_3 \pe_{C,n} \ldots,$$ it is
eventually constant in every bounded region in $\R^d$.
\label{th:constant}
\end{theorem}

Theorem~\ref{th:constant} follows from a favorable comparison between connected
saturation and $\mu$--diffuse domination.

\begin{lemma} For every $n$ and $K$, there exists a measure $\mu$ such that if
$\cA \se_{C,n} \cB$ are packings of $K$, then
$(\chi_\cA * \mu) > (\chi_\cB * \mu)$ everywhere.
\label{l:mu}
\end{lemma}
\begin{proof} Assume that $K$ has diameter 1.  Then
the maximum diameter of a connected union of $2n-1$ copies of $K$
is $2n-1$.  Let $\mu$ be the measure
with density function
$$\mu(dx) = \left(\frac{n}{n+1}\right)^{|x|/(2n-1)}dx.$$
If $B$ is the union of $k<n$ interior-disjoint copies of $K$, $A$ is the union
of $k+1$ such copies, and $A \cup B$ is connected, then we claim that
$$\mu(A) \ge \mu(B).$$
The reason is that if we express $\mu(A)$ and $\mu(B)$ as integrals,
the ratio of volumes of $A$ to $B$ is
$$\frac{k+1}{k} > \frac{n}{n+1},$$
but the ratio of the integrands is at least $n/(n+1)$.

Suppose that $\cA$ is formed from $\cB$ by replacing the components of $B$ with
the components of $A$. Then
$$(\chi_\cA * \mu)(p) - (\chi_\cB * \mu)(p) = \mu(A - p) - \mu(B - p) > 0,$$
as desired.
\end{proof}

Note that Lemma~\ref{l:mu} implies that a $\mu$--dominant packing is necessarily
$n$--saturated.  Thus a diffusively dominant packing is completely saturated. By
the argument of the lemma, if a packing $\cA$ maximizes $\mu * \chi_\cA$ at a
single point, then $\cA$ is already $n$--saturated. This produces $n$--saturated
packings which may be far from periodic.  If $\mu$ is determined by
Lemma~\ref{l:mu}, then $\mu$--dominance is slightly stronger than
$n$--saturation, but one can find such a $\mu$--dominant packing by maximizing
$\nu * \mu * \chi_cA$ at any single point, where $\nu$ is any measure with full
support on $\R^d$. (For example, $\nu = \mu$.)  Alternatively, one can find
$\mu$--dominant packings by using the fact that the relation $\se_\mu$ is closed
in the space of pairs of packings and appealing to Zorn's Lemma as in
Section~\ref{s:recur}.

\begin{proof}[Proof of Theorem~\ref{th:constant}] The argument of
Lemma~\ref{l:mu} actually shows that its conclusion is ``true by a margin'', in
the following sense: Let $R$ be a bounded domain and consider the measure $\mu$
from the lemma. Then there is a constant $\epsilon>0$ such that if $\cA$ is
obtained from $\cB$ by a connected replacement of $k<n$ copies of $K$ by $k+1$,
and if this replacement intersects $R$, then the diffused measure in $B(p,r)$
increases by at least $\epsilon$:
$$(\chi_\cA * \mu)(R) \ge (\chi_\cB * \mu)(R) + \epsilon$$
If the connected replacement does not intersect $R$, the diffused measure in
$R$ at least does not decrease. At the same time, the total measure in $R$ is
bounded above by $|\mu|(\Vol R)$.  Thus only finitely many of the replacements
in the sequence
$$\cP_1 \pe_{C,n} \cP_2 \pe_{C,n} \cP_3 \pe_{C,n} \ldots$$
may meet $R$.
\end{proof}

\begin{corollary} For every integer $n$, every Euclidean body $K$ has an
$n$--satur\-ated packing which is also uniformly dense.
\end{corollary}

\begin{proof} Apply the process of Theorem~\ref{th:constant} to a uniformly
dense packing.  The usual notion of uniform density is that the limit
$$\lim_{r \to \infty} \frac{\chi_{\cP}(B(p,r))}{\Vol B(p,r)}$$
converges uniformly in $p$.  If $\mu$ is a finite measure, then it is
equivalent to demand that
\begin{equation} |\mu|\delta(\cP) =
\lim_{r \to \infty} \frac{(\mu * \chi_{\cP})(B(p,r))}{\Vol B(p,r)}
\label{e:diffuniform}\end{equation}
converges uniformly in $p$. Since the diffused measure $\mu * \chi_{\cP}$
strictly increases under connected replacements of $k < n$ copies, and since
the density is already maximized, the limit in equation~\eqref{e:diffuniform}
does not change and continues to converge uniformly.
\end{proof}

\begin{fullfigure}{f:mass}
    {Packing a square with a mass in one corner}
\makebox[1.5em]{$\cA$} \pspicture*[](-5.25,-3.75)(5.25,3.75)
\multirput(-5.25, 3.75)(1.5,0){4}{\seslug}
\multirput(  .75, 3.75)(1.5,0){4}{\swslug}
\multirput(-5.25, 2.25)(1.5,0){4}{\seslug}
\multirput(  .75, 2.25)(1.5,0){4}{\swslug}
\multirput(-5.25,  .75)(1.5,0){4}{\seslug}
\multirput(  .75,  .75)(1.5,0){4}{\swslug}
\multirput(-5.25, -.75)(1.5,0){4}{\neslug}
\multirput(  .75, -.75)(1.5,0){4}{\nwslug}
\multirput(-5.25,-2.25)(1.5,0){4}{\neslug}
\multirput(  .75,-2.25)(1.5,0){4}{\nwslug}
\multirput(-5.25,-3.75)(1.5,0){4}{\neslug}
\multirput(  .75,-3.75)(1.5,0){4}{\nwslug}
\psframe(-5.25,-3.75)(5.25,3.75)
\endpspicture \\ \vspace{\baselineskip}
\makebox[1.5em]{$\cB$} \pspicture*[](-5.25,-3.75)(5.25,3.75)
\multirput(-5.25, 3.75)(1.5,0){8}{\seslug}
\multirput(-5.25, 2.25)(1.5,0){8}{\seslug}
\multirput(-5.25,  .75)(1.5,0){8}{\seslug}
\multirput(-5.25, -.75)(1.5,0){8}{\seslug}
\multirput(-5.25,-2.25)(1.5,0){8}{\seslug}
\multirput(-5.25,-3.75)(1.5,0){8}{\seslug}
\psframe(-5.25,-3.75)(5.25,3.75)
\endpspicture
\end{fullfigure}

\begin{example} Consider a square rim with an extra mass in one corner, as in
Figure~\ref{f:mass}.   The figure shows a pair of packings $\cA$ and $\cB$ such
that $\cB$ is uniformly recurrent, $\cA$ is not, and $\cB$ is a limit of $\cA$.
At the same time, we conjecture that $\cA$ is diffusively dominant, while $\cB$
is certainly not, since $\cA$ strictly diffusively dominates $\cB$.
\end{example}

The example of a square with a mass in one corner illustrates an artificial
aspect of diffusive dominance:  It depends on the mass distribution of the body
$K$, which is not naturally determined by the geometry of $K$.  This
shortcoming is absent for packings by translates. Let $C$ be the set of centers
of some packing of translates of some body $K$.  Let $\chi_C$ denote the
measure with a unit atom at each point in $C$.  An arbitrary mass distribution
for $K$ may be represented by a finite measure $\mu$, and the corresponding
mass distribution of the packing is $\mu * \chi_C$.  If $\nu$ is another such
measure, then $\mu * \chi_C$ is diffusively dominant if and only if $\nu *
\chi_C$ is, since they become equal if we convolve the former with $\nu$ and
the latter with $\mu$.

\begin{conjecture} The hexagonal packing is the unique diffusively dominant
circle packing.
\end{conjecture}

More generally, we conjecture that among packings by translates of a centrally
symmetric convex body in $\R^2$, the densest lattice packings are diffusively
dominant, and no other packings are.

\section{Others' notions}
\label{s:others}

In this section we relate our notions to Fejes T\'oth's solid packings
and coverings \cite{Fejes-Toth:solid} and Conway and Sloane's definition
of tight sphere packings \cite{CS:best}.

\begin{fullfigure}{f:solid}{A solid tiling which is not weakly recurrent}
\pspicture*[](-5.25,-3.75)(5.25,3.75)
\multirput(-5.25,-3)(0,1.5){3}{\psline(0,0)(4.35,0)(4.5,-.25)(4.65,0)(10.5,0)}
\multirput(-5.25,1.5)(0,1.5){2}{\psline(0,0)(4.35,0)(4.5,.25)(4.65,0)(10.5,0)}
\multirput(-4.5,-3.75)(1.5,0){3}{\psline(0,0)(0,4.35)(-.25,4.5)(0,4.65)(0,7.5)}
\multirput(0,-3.75)(1.5,0){4}{\psline(0,0)(0,4.35)(.25,4.5)(0,4.65)(0,7.5)}
\psframe(-5.25,-3.75)(5.25,3.75)
\endpspicture
\end{fullfigure}

The shape packed in Figure~\ref{f:mass} does not admit a solid packing, while
the tiling of the three shapes in Figure~\ref{f:solid} is solid but not even
weakly recurrent. Our examples in two dimensions suggest that solidity diverges
from uniform or weak recurrence for sphere packings in high dimensions.  The
relationship with diffusive dominance is less clear.  Certain sphere packings,
in particular the hexagonal circle packing, the $E_8$ lattice, and the Leech
lattice \cite{Cohn:thesis}, might have special properties that imply that they
are simultaneously solid, uniquely uniformly recurrent, and diffusively
dominant.

Conway and Sloane say the following:

\begin{quote}
Suppose we can dissect the space of a packing into finitely many polyhedral
pieces in such a way that each sphere center lies in the interior of some
piece, and there are also some \emph{empty} pieces containing no centers. Then
if we can rearrange the \emph{nonempty} pieces into another dissection in which
the centers are at least as far apart as they were originally, we call the
packing loose.

Provisionally, we may call a packing tight if it is not loose. Packings that
are tight in this sense certainly have the highest possible density. However,
we are not sure that this particular definition is the right one, and perhaps
some other meaning for ``tight'' should be used in the Postulates below.
\end{quote}

We argue that this definition is indeed not the right one. We provisionally
call the Conway--Sloane operation a \emph{tightening}.  We claim that any
periodic sphere packing in two or more dimensions can be tightened to introduce
a hole.  Figure~\ref{f:tighten} shows the construction for the hexagonal circle
packing in $\R^2$. We slide $B$ and $D$ past $A$, $C$, $E$, and $F$ to create
space for two circle centers.  Then we slide $C$ and $D$ past $A$ and $B$,
discarding $F$ and using $E$ to plug the gap between $A$ and $B$. In the end
one circle center has disappeared.

\begin{fullfigure}{f:tighten}{Tightening the hexagonal circle packing}
\pspicture*(-6,-3.464)(6,3.464)
\multirput(-5.5,2.598)(1,0){12}{\qdisk(0,0){.05}}
\multirput(-5,1.732)(1,0){11}{\qdisk(0,0){.05}}
\multirput(-5.5,.866)(1,0){12}{\qdisk(0,0){.05}}
\multirput(-5,0)(1,0){11}{\qdisk(0,0){.05}}
\multirput(-5.5,-.866)(1,0){12}{\qdisk(0,0){.05}}
\multirput(-5,-1.732)(1,0){11}{\qdisk(0,0){.05}}
\multirput(-5.5,-2.598)(1,0){12}{\qdisk(0,0){.05}}
\psframe(-6,-3.464)(6,3.464)
\psline(-6,-1.016)(-3.65,-1.016)(-3.65,-.716)(6,-.716)
\rput{-40.9}(0,1.732){\psline(1.85,.15)(1.85,-.15)
    \psline(-5,.15)(2.496,.15)(2.496,-.15)
    \psline(-.15,.15)(-.15,-.15)(3.57,-.15)}
\rput{-40.9}(-2,1.732){\psline(3.57,-.15)(10,-.15)}
\rput(-2.5,1.299){\rnode{a}{$A$}}\pnode(-3.5,1.299){a1}
\rput(-1.5,-2.165){\rnode{b}{$B$}}\pnode(-0.5,-2.165){b1}
\rput(4,1.299){\rnode{c}{$C$}}
\pnode([angle=139.1,nodesep=.6]c){c1}\pnode(3,1.299){c2}
\rput(4,-2.165){\rnode{d}{$D$}}
\pnode([angle=139.1,nodesep=.6]d){d1}\pnode(5,-2.165){d2}
\rput(-.5,1.299){\rnode{e}{$E$}}\pnode(.5,1.299){e1}
\rput(2.667,.289){\rnode{f}{$F$}}\pnode(1.667,.289){f1}
\ncline[nodesepA=.15]{->}{b}{b1}
\ncline{->}{c}{c1}
\ncline{->}{d}{d1}\ncline[nodesepA=.15]{->}{d}{d2}
\ncline[nodesepA=.15]{e}{e1}
\ncline[nodesepA=.15]{f}{f1}
\endpspicture\end{fullfigure}

We can generalize the argument of Theorem~\ref{th:hex} to produce a sensible
criterion which is similar to the Conway--Sloane proposal.  Let $K$ be a domain
with unit volume and with packing density $\delta$.  Suppose that for each
packing $\cP$ of $K$, we can form a packing $\cQ$ of domains (in general
non-congruent) with volume $1/\delta$ whose elements are in bijection with
those of $\cP$.  Suppose further that $\cQ$ satisfies the following conditions.

\begin{description}
\item[\rm(1)] The assignment $\cP \mapsto \cQ$ is translation-invariant
on the space of packings.
\item[\rm(2)] The assignment $\cP \mapsto \cQ$ is continuous in the
Hausdorff topology on the space of packings.
\item[\rm(3)] Each element of $\cQ$ lies a bounded distance from the associated
element of $\cP$, with the bound independent of $\cP$.
\end{description}

Then:

\begin{theorem} If $\cP$ and $\cQ$ are defined as above, and if
$\cP$ is weakly recurrent and dense, then $\cQ$ is a tiling.
\label{th:nogap}
\end{theorem}

The complicated hypothesis of Theorem~\ref{th:nogap} is designed to fit the
argument of Theorem~\ref{th:hex}.  Theorem~\ref{th:hex} is a special case
because we can form the packing $\cQ$ from the Voronoi tiling of a circle
packing.  For example, we can truncate each Voronoi region so that it has the same
area as the circumscribed hexagon.

In conclusion, Theorem~\ref{th:nogap} provides a sense, as Conway and Sloane
desired, in which weakly recurrent dense packings have no gaps
\cite{Sloane:confirmed}.  So we may adopt their postulate as a conjecture:

\begin{conjecture} Every weakly recurrent, dense sphere packing in dimension $2
\le n \le 8$ fibers over one in dimension $2^k$, where $2^k$ is the largest
power of 2 strictly less than $n$.
\end{conjecture}

\end{document}